\documentclass[12pt, a4paper]{article}
\usepackage[latin1]{inputenc}
\usepackage[english]{babel}
\usepackage{amsmath}
\usepackage[english]{babel}
\usepackage{indentfirst}
\usepackage{amstext}
\usepackage{amsfonts}
\usepackage{textcomp}
\usepackage{amssymb}
\usepackage[dvips]{graphicx}
\usepackage{setspace}
\usepackage[all]{xy}

\newcommand {\demo}{\hskip -0.6cm{\bf Proof.  }}
\newcommand {\fim}{\hfill{$\square$}\vskip 1pc}

\newcommand {\N}{\mathbb{N}}
\newcommand {\Z}{\mathbb{Z}}
\newcommand{\s}{\sigma}

\newcommand {\F}{\mathbb{F}}

\newtheorem{teorema}{Theorem}[section]
\newtheorem{lema}[teorema]{Lemma}

\newtheorem{definicao}[teorema]{Definition}
\newtheorem{proposicao}[teorema]{Proposition}
\newtheorem{exemplo}[teorema]{Example}
\newtheorem{remark}[teorema]{Remark}
\begin{document}
\onehalfspace

\title{Ultragraphs and shift spaces over infinite alphabets}
\maketitle
\begin{center}

{\large Daniel Gon\c{c}alves\footnote{This author is partially supported by CNPq and Capes PVE085/2012.} and Danilo Royer}\\
\end{center}  
\vspace{8mm}

\abstract 
In this paper we further develop the theory of one sided shift spaces over infinite alphabets, characterizing one-step shifts as edge shifts of ultragraphs and partially answering a conjecture regarding shifts of finite type (we show that there exists shifts of finite type that are not conjugate, via a conjugacy that is eventually finite periodic, to an edge shift of a graph ).  We also show that there exists edge shifts of ultragraphs that are shifts of finite type, but are not conjugate to a full shift, a result that is not true for edge shifts of graphs. One of the key results needed in the proofs of our conclusions is the realization of a class of ultragraph C*-algebras as partial crossed products, a result of interest on its own. 
\doublespace\newline\newline
MSC 2010: 37B10, 54H20 
\doublespace\newline
Keywords: Symbolic dynamics, one-sided shift spaces, infinite alphabets, ultragraph C*-algebras, partial crossed products.

\section{Introduction}

Symbolic dynamics is one of the most important areas in dynamical systems and is usually concerned with a finite alphabet $A$, spaces of infinite sequences $A^\N$ or $A^\Z$ in $A$, and subspaces of $A^\N$ or $A^\Z$ which are closed under the shift map. Over the last decades there has been many generalizations proposed for the infinite alphabet case (see \cite{OTW} for an overview). Recently Ott, Tomforde and Willis, see \cite{OTW}, proposed a generalization $\Sigma_A$ of $A^\N$ which can be identified with the set of all finite sequences in $A$ (including an empty sequence) together with the infinite sequences in $A$. Subshifts in this case are then taken to be closed subsets $\Lambda\subset \Sigma_A$ invariant under the shift which satisfy the so called ``infinite extension property'', which guarantees that $\Lambda$ contains infinite sequences and is determined by its infinite sequences. 

The approach in \cite{OTW} for one sides shifts over infinite alphabets (from now on called OTW approach) was taken with an eye towards the interplay with C*-algebras of infinite graphs and a comprehensive study of such shifts and morphisms between them is presented in \cite{OTW}, along with two conjectures  left open. One of the conjectures regarded the existence of $M+1$ step shifts that are not conjugate to an $M$ step shift, which we proved to be true in \cite{Mstep}, and the second one, namely the existence of shifts of finite type that are not conjugate to an edge shift of a graph is the motivating question behind this paper. It is interesting to note that recently there has been a lot of activity regarding the OTW approach for shifts over infinite alphabets. In particular Curtis-Hedlund-Lyndon type theorems were proved in \cite{GSS} and \cite{GS}, two sided shift spaces over infinite alphabets were defined in \cite{GSS1} and semigroup structures over OTW shifts were studied in \cite{GSS2}. 

In the course of our attempts to prove the conjecture regarding shifts of finite type we realized that the theory of ultragraphs (ultragraphs are generalizations of graphs introduced by Mark Tomforde, see \cite{Tom}) would be helpful, as it would provide a link between dynamical systems and C*-algebras, hopefully working as the link between one sided shift spaces and graph algebras explored in \cite{OTW}. So we introduce edge shifts associated to ultragraphs in Section 3, where we also show that the class of ultragraph edge shifts coincide with the class of 1-step OTW shift spaces. Furthermore, assuming that there exists a conjugacy between the edge shift of the graph $E$ and a shift of finite type, we give a partial description of $E$ and, under the additional hypothesis the the conjugacy is eventually finite periodic preserving (a mild condition we define in Section 3), we show that the graph $E$ is a bouquet of loops and therefore the associated edge shift is the full shift. This is as far as we are able to proceed using only topological dynamics arguments. 

To further develop the solution of our problem we envisioned the use of the interplay between Dynamical Systems and C*-algebras. The main idea is to prove that if two edge shifts from ultragraphs are conjugate then their ultragraph C*-algebras are isomorphic and hence so are their K-theory. So if a SFT (given as the edge shift of an ultragraph) is conjugate to the edge shift of a graph, then the ultragraph C*-algebras associated should be isomorphic and hence have the same K-theory. We alert the reader that we were not able to prove that ultragraph C*-algebras are invariant under conjugacy of ultragraph edge shifts in full generality, but instead we were able to prove it for a class large enough for our purposes regarding the SFT conjecture. Furthermore, to prove the invariance for this class of ultragraphs, we realize these ultragraph C*-algebras as partial crossed products, a result we believe to be of interest on it own.

Now, by the results from section 3, we know that the edge shift of a graph that is conjugate, via a eventually finite periodic preserving conjugacy, to a shift of finite type must be a full shift. Hence, to prove the SFT conjecture (for finite periodic preserving conjugacies), all we need is to find an edge ultragraph SFT whose associated ultragraph C*-algebra has K-theory different from the K-theory of the ultragraph (in this case graph) C*-algebra associated to the full shift. 

We notice that if conjugacy preserves shifts of finite type, that is, the image of a shift of finite type by a conjugacy is also a shift of finite type, then we can drop the eventually finite periodic preserving condition in our results, as the theory developed in section 3 would then readily imply that an edge shift of a graph conjugate to a shift of finite type must be the full shift. Unfortunately we do not know the answer to this question and leave it open: If $X$ is a shift of finite type and $\phi:X \rightarrow Y$ is a conjugacy, is $Y$ a shift of finite type?

A final interesting observation we derive from the consideration of ultragraph C*-algebras is that, differently from edge shifts from graphs, there exists ultragraphs edge shifts of finite type that are not conjugate to a full shift.

Below we describe the organization of the paper:

As we already mentioned, in section 3, we study ultragraph edge shifts. In section 4, building from ideas in \cite{GR1} for Leavitt path algebras, we realize a class of ultragraph C*-algebras as partial crossed products C*-algebras and, using partial actions theory, show that, for a class of ultragraphs, conjugacy between ultragraph edge shifts implies isomorphism of the ultragraph C*-algebras associated. Finally, in section 5, we present two examples of shifts of finite type, given as the edge shift of an ultragraph, such that the K-theory of the ultragraph C*-algebras associated has torsion and hence it is different from the K-theory of the ultragraph C*-algebra associated to the full shift, therefore proving the existence of shifts of finite type not conjugate, via a eventually finite periodic preserving conjugacy, to an edge shift of a graph. Since the examples we construct are shifts of finite type, the arguments in this section also prove that there exists ultragraphs edge shifts of finite type that are not conjugate to a full shift. To make the paper as self contained as possible we recall the OTW approach and ultragraph C*-algebras in a background section, which we present below.

\section{Background}

\subsection{One sided shift spaces over infinite alphabets}\label{OTWapproach}

In this subsection we recall the OTW approach to one sided shift spaces over infinite alphabets as in \cite{OTW}.

Given a set (alphabet) $A$, define a new symbol, say $\O$, which will indicate the {\em empty sequence} and let
$\Sigma_A^{fin}:=\{\O\}\cup\bigcup\limits_{k\in\N}A^k$ and $\Sigma_A^{inf}:=\{(x_i)_{i\in\N}:\ x_i\in A\ \forall i\in\N\}$.
We say that $\Sigma_A^{fin}$ is the set of all finite sequences over $A$, while $\Sigma_A^{inf}$ is the set of all infinite sequences over $A$.

The {\em full shift} over $A$ is the set $$\Sigma_A:=\left\{\begin{array}{lcl}\Sigma_A^{inf} & \text{ if} & |A|<\infty\\ \Sigma_A^{inf}\cup\Sigma_A^{fin} & \text{ if} & |A|=\infty.\end{array}\right.$$
Given $x\in\Sigma_A$ we define $l(x)$ (also denoted by $|x|$) as the length of $x$, which is equal to $k$ if $x\in A^k$ and equal to $\infty$ otherwise.

Next we recall the basis for the topology in $\Sigma_A$. Given $x=(x_1,\ldots,x_k)\in\Sigma_A^{fin}$, $x\neq\O$, and a finite set $F\subset A$, we define a {\em generalized cylinder set} of $\Sigma_A$ as
$Z(x,F):=\{y\in\Sigma_A:\ y_i=x_i\ \forall i=1\ldots,k,\ y_{k+1}\notin F\}$. If $ x = \O $ then we let $Z(\O,F)= \{y\in\Sigma_A:\ y_1\notin F\}$.
														
If $F=\emptyset$ we will denote $Z(x,F)$ just by $Z(x)$. Notice that $Z(x)$ coincides with the usual cylinders which form the basis of the product topology in $\Sigma_A^{inf}$, while $Z(\O,F)$ is a neighborhood of the empty sequence $\O$. It is proven in \cite{OTW} that the collection of generalized cylinder sets defined above forms a basis for a topology in $\Sigma_A$ consisting of compact open sets. The full shift space over $A$ is defined as $\Sigma_A$ with the topology above and in \cite{OTW} it is proved that $\Sigma_A$ is zero dimensional (that is, it has a basis of clopen sets) and compact.

The {\em shift map} $\s:\Sigma_A\to\Sigma_A$ is defined as follows:
$$\s(x)=\left\{\begin{array}{lcl} \O & \text{ if} & x=\O \text{ or } x\in A^1\\
                      (x_{i+1})_{i=1}^{k-1} & \text{ if} & x=(x_i)_{i=1}^k\in A^k\\
                 (x_{i+1})_{i\in\N} & \text{ if} & x=(x_i)_{i\in\N}\in \Sigma_A^{inf}.\end{array}\right.$$ 


We will be particularly interested in special subsets of $\Sigma_A$, which are named shift spaces.
We say that $\Lambda\subseteq\Sigma_A$ is a {\em shift space} over $A$ if the following three properties hold:
\begin{enumerate}
\item[1] $\Lambda$ is closed with respect to the topology of $\Sigma_A$;
\item[2] $\Lambda$ is invariant under the shift map, that is, $\s(\Lambda)\subseteq\Lambda$;
\item[3] $\Lambda$ satisfies the `infinite-extension property', that is, for all $x \in \Lambda$, with $l(x)< \infty$, the set $\{a \in A: xay \in \Lambda \text{ for some } y \in \Sigma_A\}$ is infinite.
\end{enumerate}

Notice that, by Proposition~3.8 of \cite{OTW}, $\Lambda^{inf}$ is dense in $\Lambda$.

In symbolic dynamics, an equivalent way to define a shift space is via a set of forbidden words: Let $\mathbf{F}\subset \bigcup\limits_{k\geq 1}A^k$, 

$\begin{array}{ll} X_{\mathbf{F}}^{inf} = & \{ x \in \Sigma_A^{inf}: \text{ no subblock of $x$ is in $\mathbf{F}$} \}, \text{ and } \\ 
X_{\mathbf{F}}^{fin} = & \{ x \in \Sigma_A^{fin}: \text{ there are infinitely many } a \in A \\
& \text{ for which there exists } y\in \Sigma_A^{inf} \text{  such that } xay\in  X_{\mathbf{F}}^{inf}  \},
\end{array}$
where a subblock of $x$ means an element $u\in \Sigma_A^{fin}$ such that $x=yuz$ for some $y\in \Sigma_A^{fin}$ and some $z \in \Sigma_A$. Now, define $ X_{\mathbf{F}}$ as $X_{\mathbf{F}}^{inf} \cup X_{\mathbf{F}}^{fin}$. It is proven in \cite{OTW}, theorem 3.16, that $\Lambda$ is a shift space if and only if $\Lambda = X_{\mathbf{F}}$ for some $\mathbf{F}\subset \bigcup\limits_{k\geq 1}A^k$.

Following \cite{OTW}, we say that $\Lambda$ is a shift of finite type, SFT, if we can take $\mathbf{F}$ having only finitely elements.


Finally recall that a map $\phi: \Lambda \rightarrow Y$, where $\Lambda\subseteq \Sigma_A$, $Y\subseteq \Sigma_B$ are shift spaces over alphabets $A$ and $B$ respectively, is a conjugacy if, and only if, it is bijective, continuous, shift commuting and $l(\phi(x)) = l(x)$ for all $x \in \Lambda$. 

\subsection{Ultragraph C*-algebras}

Ultragraph C*-algebras were introduced by Tomforde in \cite{Tom}. Here we recall the main definitions and relevant results.

\begin{definicao}\label{def of ultragraph}
An \emph{ultragraph} is a quadruple $\mathcal{G}=(G^0, \mathcal{G}^1, r,s)$ consisting of two countable sets $G^0, \mathcal{G}^1$, a map $s:\mathcal{G}^1 \to G^0$, and a map $r:\mathcal{G}^1 \to P(G^0)\setminus \{\emptyset\}$, where $P(G^0)$ stands for the power set of $G^0$.
\end{definicao}

\begin{exemplo}\label{exemplo1} Let $\mathcal{G}$ be the ultragraph with a countable number of vertices, say $G^0=\{v_i\}_{i\in \N}$, and edges such that $s(e_i)= v_i$ for all $i$, $r(e_1)=\{v_3, v_4, v_5, \ldots \}$ and $r(e_j)= G^0$ for all $j\neq 1$. We can represent this ultragraph as in the picture below. 
\end{exemplo}

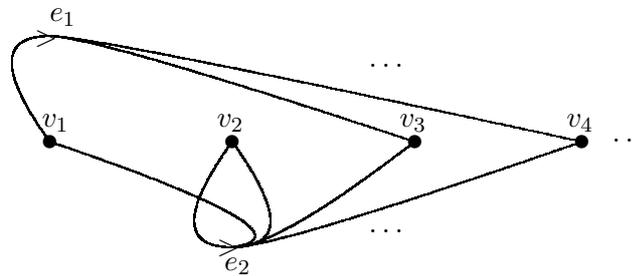
\centerline{
\setlength{\unitlength}{2cm}
\begin{picture}(4,1.2)
\put(0.5,0){\circle*{0.08}}
\put(0.5,0){\qbezier(0,0)(-0.5,0.7)(0,0.7)}
\put(0.4,0.64){$>$}
\put(0.5,0.7){\qbezier(0,0)(0.3,0)(2.4,-0.7)}
\put(0.5,0.7){\qbezier(0,0)(0.3,0)(3.5,-0.7)}
\put(1.7,0){\qbezier(0,0)(-0.5,-0.7)(0,-0.7)}
\put(1.7,0){\qbezier(0,0)(0.5,-0.7)(0,-0.7)}
\put(1.6,-0.75){$>$}
\put(1.7,-0.7){\qbezier(0,0)(0.3,0)(1.2,0.7)}
\put(1.7,-0.7){\qbezier(0,0)(0.3,0)(2.3,0.7)}
\put(1.7,-0.7){\qbezier(0,0)(0.6,0.085)(-1.2,0.7)}
\put(0.45,0.1){$v_{1}$}
\put(1.7,0){\circle*{0.08}}
\put(1.6,0.1){$v_2$}
\put(2.9,0){\circle*{0.08}}
\put(2.8,0.1){$v_3$}
\put(4,0){\circle*{0.08}}
\put(3.9,0.1){$v_4$}
\put(4.2,0){\dots}
\put(2.6,0.5){\dots}
\put(2.6,-0.6){\dots}
\put(0.5,0.8){$e_1$}
\put(1.65,-0.85){$e_2$}
\end{picture}}
\vspace{2 cm}

Before we define the C*-algebra associated to an ultragraph we need the following notion.

\begin{definicao}\label{def of mathcal{G}^0}
Let $\mathcal{G}$ be an ultragraph. Define $\mathcal{G}^0$ to be the smallest subset of $P(G^0)$ that contains $\{v\}$ for all $v\in G^0$, contains $r(e)$ for all $e\in \mathcal{G}^1$, and is closed under finite unions and nonempty finite intersections.
\end{definicao}

\begin{definicao}\label{def of C^*(mathcal{G})}
Let $\mathcal{G}$ be an ultragraph. The \emph{ultragraph algebra} $C^*(\mathcal{G})$ is the universal $C^*$-algebra generated by a family of partial isometries with orthogonal ranges $\{s_e:e\in \mathcal{G}^1\}$ and a family of projections $\{p_A:A\in \mathcal{G}^0\}$ satisfying
\begin{enumerate}
\item\label{p_Ap_B=p_{A cap B}}  $p_\emptyset=0,  p_Ap_B=p_{A\cap B},  p_{A\cup B}=p_A+p_B-p_{A\cap B}$, for all $A,B\in \mathcal{G}^0$;
\item\label{s_e^*s_e=p_{r(e)}}$s_e^*s_e=p_{r(e)}$, for all $e\in \mathcal{G}^1$;
\item $s_es_e^*\leq p_{s(e)}$ for all $e\in \mathcal{G}^1$; and
\item\label{CK-condition} $p_v=\sum\limits_{s(e)=v}s_es_e^*$ whenever $0<\vert s^{-1}(v)\vert< \infty$.
\end{enumerate}
\end{definicao}

\begin{remark}
The ultragraphs we use to prove theorem \ref{provafracaconjectura} satisfy condition (L), that is, every loop has an exit (a loop $\alpha:=\alpha_1 \cdots \alpha_n$ has an exit if there is either  an edge $e\in \mathcal{G}^1$ such that $s(e)\in r(\alpha_i)$ and $e\neq \alpha_{i+1}$ for some i, or there is a sink $w$ with $w \in r(\alpha_i)$ for some $i$). Under this condition the Cuntz-Krieger uniqueness theorem, which we state below, holds. 
\end{remark}


\begin{teorema}\cite[Theorem~6.7]{Tom} Let  $\mathcal{G}$  be an ultragraph that satisfies condition (L) and let $\pi$ be a representation of $C^*(\mathcal{G})$ such that $\pi(P_A) \neq 0$ for all $A \in  \mathcal{G}^0$. Then $\pi$ is faithful.
\end{teorema}

\section{The edge shift of an ultragraph}

In \cite{OTW} a comprehensive study of shifts of one sided shift spaces over infinite alphabets is presented. In this section we further develop the theory, by introducing ultragraph edge shift spaces. 

Let $\mathcal{G}=(G^0,\mathcal{G}^1, r,s)$ be an ultragraph. Consider the set of edges as the alphabet, that is, let $A=\mathcal{G}^1$ and let $\Sigma_A$ denote the full shift as in section \ref{OTWapproach}.

\begin{definicao}Let $\mathcal{G}=(G^0,\mathcal{G}^1, r,s)$ be an ultragraph and let $\mathcal{G}^\infty=\{(x_n)_{n\in \N}\in A^\N: s(x_{i+1})\in r(x_i) \ \ \forall i\in \N\}$, that is, $\mathcal{G}^\infty$ is the set of all infinite paths in $\mathcal{G}$. The edge shift associated to the ultragraph $\mathcal{G}$ (also called the ultragraph edge shift of $\mathcal{G}$) is the closure of $\mathcal{G}^{\infty}$ in $\Sigma_A$.
\end{definicao}

\begin{remark}
When $\mathcal{G}$ is a graph then the definition above reduces to the definition of an edge shift given in \cite{OTW}.
\end{remark}

It was shown in \cite{OTW} that each edge shift of a graph is a 1-step shift, but the converse does not hold (see \cite[5.16]{OTW} and \cite[5.18]{OTW}). With the notion of ultragraph edge shifts we can completely characterize 1-step shift spaces:

\begin{proposicao} Each ultragraph edge shift is a 1-step shift and each 1-step shift is an ultragraph edge shift.
\end{proposicao} 

\demo Let $Y=\overline{\mathcal{G}^{\infty}}$ be the edge shift associated to an ultragraph $\mathcal{G}$. We will show that $Y$ a 1-step shift space. For this let $F\subseteq \Sigma_{\mathcal{G}^1}$ be defined by $$F=\{ef:e,f\in \mathcal{G}^1, s(f)\notin r(e)\}.$$ Recall that $X_F=X_F^{inf}\cup X_F^{fin}$, where  
$$X_F^{inf}=\{x\in \Sigma_{\mathcal{G}^1}^{inf}:\text{ no subblock of $x$ is in $F$ }\}$$ and 
$$X_F^{fin}=\{x\in \Sigma_{\mathcal{G}^1}^{fin}: \text{there are infinitely many } a\in \mathcal{G}^1$$ 
$$\hspace{2cm}\text{such that $xay\in X_F^{inf}$ for some $y\in \Sigma_{\mathcal{G}^1}^{inf}$} \}.$$
By \cite[3.14]{OTW}, $X_F$ is a 1-step shift space. Furthermore, by \cite[3.8]{OTW}, $X_F^{inf}$ is dense in $X_F$ and, since $X_F$ is closed, we have that $\overline{X_F^{inf}}=X_F$. Now notice that $\mathcal{G}^{\infty}=X_F^{inf}$, so that $Y=\overline{\mathcal{G}^{\infty}}=\overline{X_F^{\inf}}=X_F$ and hence $Y$ is a 1-step shift space.

For the converse, consider an alphabet $A$ and let $X_F$ be a 1-step shift space, where $F\subseteq \Sigma_A$ is a set of words of length $2$. Let $\mathcal{G}=(G^0,\mathcal{G}^1, r,s)$ be the ultragraph defined by $G^0=\{v_a:a\in A\}$, $\mathcal{G}^1=A$, and for all $a\in A$, $s(a)=v_a$ and $r(a)=\{v_b: ab\notin F\}$. By definition $\mathcal{G}^{\infty}=\{(x_n)_{n\in \N}\in A^\N: s(x_{i+1})\in r(x_i)\}$ and by the first part of this proof $\overline{\mathcal{G}^{\infty}}$ is in fact a shift space. Furthermore, $$\mathcal{G}^{\infty}=\{(x_n)_{n\in \N}\in A^\N: s(x_{i+1})\in r(x_i)\,\,\,\forall i\in \N\}=$$ 
$$=\{(x_n)_{n\in \N}\in A^\N: v_{x_{i+1}}\in \{v_b:x_ib\notin F\}\forall i\in \N\}=$$ $$=\{(x_n)_{n\in \N}\in A^\N:x_ix_{i+1}\notin F\,\,\,\forall i\in \N\}=X_F^{\inf},$$ and it follows from \cite[3.9]{OTW} that $\overline{\mathcal{G}^{\infty}}=X_F.$
\fim

Since each graph $E$ is an ultragraph, each edge shift of a graph is an edge shift of an ultragraph. Moreover, by the previous proposition and by \cite[5.18]{OTW}, there are edge shits of ultragraphs which are not edge shifts of graphs and so we get that
$$\{\text{ edge shifts of graphs }\}\subsetneq\{\text{ edge shifts of ultragraphs }\}.$$

Next we further characterize ultragraph edge shifts. For this let, for each $n\in \N$, $n>0$, $\mathcal{G}^n$ be the set of all paths of length $n$, that is, the set of all the words $\alpha_1....\alpha_n$ such that $s(\alpha_{i+1})\in r(\alpha_i)$ for each $i\in \{1,...,n-1\}$.

\begin{proposicao}\label{p1} Let $\mathcal{G}$ be an ultragraph and $X_\mathcal{G}$ be the associated edge shift. Suppose that for each $e\in \mathcal{G}^1$ there exists a vertex in $r(e)$ which is not a sink. If $\mathcal{G}^1$ is finite then $X_{\mathcal{G}}=\mathcal{G}^{\infty}$ and if $\mathcal{G}^1$ is infinite then $$X_{\mathcal{G}}=\mathcal{G}^{\infty}\cup{\{\O\}}\cup\{\alpha \in \bigcup\limits_{n=1}^\infty\mathcal{G}^n: s^{-1}(r(\alpha)) \text{ is infinite}\}.$$
\end{proposicao}

\demo By definition, $X_\mathcal{G}=\overline{\mathcal{G}^\infty}$. If $\mathcal{G}^1$ is finite then $\mathcal{G}^\infty$ is closed, and so $X_\mathcal{G}=\mathcal{G}^\infty$. Now suppose $\mathcal{G}^1$ is infinite. First we show that 
$$X_{\mathcal{G}}\supseteq\mathcal{G}^{\infty}\cup{\{\O\}}\cup\{\alpha\in \bigcup\limits_{n=1}^{\infty}\mathcal{G}^n: s^{-1}(r(\alpha)) \text{ is infinite}\}.$$

Note that $\O\in X_\mathcal{G}$ because there are infinite edges in $\mathcal{G}^1$ and, for each edge $e\in \mathcal{G}^1$, there exists an element $(x_n)_{n\in \N}\in \mathcal{G}^\infty$ such that $x_1=e$ (since for each $f\in \mathcal{G}^1$, $r(f)$ contains some vertex which is not a sink). Moreover, for $\alpha \in \mathcal{G}^n$ such that 
$r(\alpha)$ contains sources of infinitely many edges, let $\{e_1,e_2,e_3,...\}$ be a set of such edges and, for each edge $e_n$, let $x^n\in \mathcal{G}^{\infty}$ be such that $x^n=(\alpha_1,...,\alpha_{|\alpha|},e_n,...)$. Then the sequence $(x^n)_{n\in \N}$ converges to $\alpha$.

Next we show that $$X_{\mathcal{G}}\subseteq\mathcal{G}^{\infty}\cup{\{\O\}}\cup\{\alpha\in \bigcup\limits_{n=1}^{\infty}(\mathcal{G}^1)^n: s^{-1}(r(\alpha)) \text{ is infinite}\}.$$

Let $(x^n)_{n\in \N}\subseteq \mathcal{G}^\infty$ be a sequence converging to $x\in \Sigma_{\mathcal{G}^1}$. If $|x|=\infty$ then, as usual, $x\in \mathcal{G}^{\infty}$. Suppose that $0<|x|<\infty$ and $s^{-1}(r(x))$ is finite. Let $F=s^{-1}(r(x))$. Then the generalized cylinder $Z(x,F)$ contains some $x^k$ (since $(x^n)_{n\in \N}$ converges to $x$) and so $x_1^k...x_{|x|}^k=x_1...x_{|x|}$ and $x_{|x|+1}^k\notin F$, which is impossible, since $s(x_{|x|+1}^k)\in r(x_{|x|}^k)=r(x)$. \fim

\subsection{Conjugacy between shifts of finite type and edge shifts.}

Next we give some characterizations of a graph such that the associate edge shift is conjugate to a shift of finite type. 

\begin{lema}\label{fullshift} Let $X$ be a shift of finite type over an infinite alphabet $A$. Suppose that $X$ is conjugated to an edge shift $Y_E$, of a graph $E$. Then all elements of length one in $Y_E$, seen as edges in the graph, have the same range. Furthermore, except for possibly a finite number of elements, length one elements in $Y_E$ are loops with the same range. 
\end{lema}

\demo Let $\varphi:X\rightarrow Y_E$ be a conjugacy. First we show that, for each $e,f$ of length one in $Y_E$, we have that $r(e)=r(f)$ in $E^1$.

Since $\varphi$ is a conjugacy, for each $e, f$ of length one in $Y$ there are $a,b\in A$ such that $\varphi(a)=e$ and $\varphi(b)=f$.
Since $X$ is of finite type, we can chose infinite distinct elements $c_n\in A$, $n\in \N$,  such that $x_n:=ac_nc_nc_n...\in X$ and $y_n:=bc_nc_nc_n...\in X$, for each $n\in \N$. Then the sequences $(x_n)$ and $(y_n)$ converge to $a$ and $b$, respectively. Now, notice that for each $n\in \N$, $\varphi(x_n)=\varphi(ac_nc_nc_n...)=a^n\varphi(c_nc_nc_n...)$, where $a^n\in E^1$, and, since $\varphi(x_n)$ converges to $\varphi(a)=e$, we have that $a^n=e$ for each $n$ greater than some $n_0$. Hence $\varphi(ac_nc_nc_n....)=e\varphi(c_nc_nc_n...)$ for $n\geq n_0$. Analogously, there exist $n_1\in \N$ such that $\varphi(bc_nc_nc_n...)=f\varphi(c_nc_nc_n...)$ for $n\geq n_1$. In particular, for $m\geq \max \{n_0,n_1\}$, since $e\varphi(c_mc_mc_m...)\in Y_E$ and $f\varphi(c_mc_mc_m...)\in Y_E$, we have that  $r(e)=s(\varphi(c_mc_mc_m...))=r(f)$.

Next we show that almost all elements (except maybe a finite number of then) of length one in $Y_E$ are loops with the same range. 

Let $\{a_i\}$ be an enumeration of the alphabet $A$, and let $e_i = \varphi(a_i)$. Since $\varphi$ is a conjugacy every element of length one in $Y_E$ is equal to an $e_i$ for some i. Since $X$ is a SFT, there exists $N$ such that for all $n>N$ the element $a_1a_n$ belongs to $X$. Notice that the sequence $(a_1 a_n )$ converges to $a_1$. Therefore, since $\varphi$ is a conjugacy, there exists $M>0$ such that, for all $j>M$, $$\varphi(a_1 a_j) = \varphi (a_1) \varphi (a_j) = e_1 e_j.$$ Hence, for all $j>M$, $r(e_1)=s(e_j)$ and, since we already have that $r(e_1)=r(e_j)$ for all $j$, we obtain that $e_j$ is a loop for every $j>M$. \fim

To further characterize a graph whose associated edge shift is conjugate to an SFT we need an extra condition on the conjugacy. We define this now:

\begin{definicao} Let $X$ be a shift of finite type over an infinite alphabet, $Y$ a shift space and $\phi:X\rightarrow Y$ a conjugacy. We say that $\phi$ is eventually finite periodic if there exists $p\geq 2$ such that, for all periodic element of period $p$, say $x= x_1\ldots x_px_1\ldots x_p\ldots$, we have that $\phi (x_1\ldots x_p) = \phi(x)_1 \ldots \phi(x)_p$ (where $\phi(x)_i$ denotes the $i$ entry of $\phi(x)$). We call $p$ the period of $\phi$.
\end{definicao}

\begin{proposicao}\label{finiteperiodicfullshift} Let $X$ be a shift of finite type over an infinite alphabet $A$. Suppose that $X$ is conjugated, via an eventually finite periodic conjugacy, to an edge shift $Y_E$, of a graph $E$. Then $E$ is a bouquet of loops and $Y_E$ is the full shift.
\end{proposicao}

\demo Let $\{a_i\}_{i\in \N}$ be an enumeration of the alphabet of $X$ and $\phi: X \rightarrow Y$ be an eventually finite periodic conjugacy, with period $p$. Then each element of length one in $Y_E$ can be written as $e_i := \phi(a_i)$. By proposition \ref{fullshift} we have that $r(e_i)= v$, for all i, and there exists $N_0$ such that, for all $i\geq N_0$, the edge $e_i$ is a loop with range $v$.

Suppose that there exists $e_K$ that is not a loop. By proposition \ref{fullshift} we have that $r(e_K)=v$. Since $X$ is a SFT, there exists $N>0$ such that, for all $n>N$ the element $a_K a_n a_K$ is in $X$. Since the sequence $(a_K a_n a_K)$ converges to $a_K$, and $\phi$ is a conjugacy, there exists $M>N>0$ such that, for all $n\geq M$, $\phi(a_K a_n a_K) = e_K b_n e_K$, where $b_n$ is an edge in $E$ such that $s(b_n)=r(e_K)=v$ and $r(b_n)=s(e_K)$.

Now, consider the element of period $p$ in $Y_E$ given by $$y=e_K e_{N_0}^n b_M e_K e_{N_0}^n b_M \ldots ,$$ where $e_{N_0}^n$ denotes $n$ times the edge  $e_{N_0}$ and $n$ is such that $n+2=p$.

Let $x=a_1a_2\ldots a_p a_1 \ldots a_p \ldots$ be such that $\phi(x)=y$. Since $\phi$ has period $p$ we have that $\phi(a_1 \ldots a_p) = e_K e_{N_0}^n b_M$. This implies that $ \phi (a_p) = b_M$ and hence $s(e_K)=r(b_M) =  v = r(e_K)$, a contradiction. Therefore every $e_i$ is a loop with range $v$.

Finally, if $(e_i)$ is a sequence in $Y_E$ and one of the $e_i$, say $e_K$, belongs to $Y_E$, then $e_K$ is a loop with source $v$. Hence $r(e_{K-1})$ is also $v$, what implies that $e_{K-1}$ is also a loop. Proceeding in this fashion we see that every $e_i$, for $i\leq K$, is a loop with range $v$. Now, if none of the $e_i$ belongs to $Y_E$, let $x=(a_i)$ in $X$ be such that $\phi(x)=(e_i)$. Then $x$ can be approached by a sequence of elements of finite length, for example, $a_1$, $a_1 a_2$, $a_1 a_2 a_3$, \ldots, and hence the image of this sequence by $\phi$ approaches $y$, a contradiction, since $y$ can not be approached only by elements of finite length. 
\fim


\section{Realization of a class of ultragraph C*-algebras as partial crossed products}

Let $\mathcal{G}=(G^0,\mathcal{G}_1,r,s)$ be an ultragraph with no sinks such that $\mathcal{G}^1$ is infinite, $s^{-1}(r(e))^C=\mathcal{G}^1\setminus s^{-1}(r(e))$ is finite for each $e\in \mathcal{G}^1$ and, for each $v\in G^0$, $s^{-1}(v)$ is finite or $s^{-1}(v)^C$ is finite. In this section, we realize the C*-algebra associated to $\mathcal{G}$ as a partial crossed product (see \cite{exel} for a detailed introduction to partial actions and their associated crossed products).

Though our results hold only for ultragraph satisfying the above conditions, we can associate a partial action to any ultragraph such that  $\mathcal{G}^1$ is infinite and $s^{-1}(r(e))^C$ is finite for each $e\in \mathcal{G}^1$. So, from now on we assume that all ultragraphs satisfy these last two conditions.


\begin{exemplo}\label{ex1} Let $\mathcal{G}^1=\{e_i\}_{i\in \N}$ and $G^0=\{v_i\}_{i\in \N}$. Define $s(e_i)=v_i$ for each $i$ and $r(e_i)=\{v_j:j\ge i\}$. We show a picture of the ultragraph $\mathcal{G}$ below. 

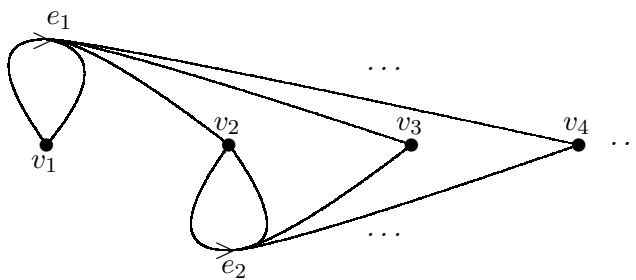
\centerline{
\setlength{\unitlength}{2cm}
\begin{picture}(4,1.2)
\put(0.5,0){\circle*{0.08}}
\put(0.5,0){\qbezier(0,0)(-0.5,0.7)(0,0.7)}
\put(0.5,0.7){\qbezier(0,0)(0.5,-0.1)(0,-0.7)}
\put(0.4,0.64){$>$}
\put(0.5,0.7){\qbezier(0,0)(0.3,0)(1.2,-0.7)}
\put(0.5,0.7){\qbezier(0,0)(0.3,0)(2.4,-0.7)}
\put(0.5,0.7){\qbezier(0,0)(0.3,0)(3.5,-0.7)}
\put(1.7,0){\qbezier(0,0)(-0.5,-0.7)(0,-0.7)}
\put(1.7,0){\qbezier(0,0)(0.5,-0.7)(0,-0.7)}
\put(1.6,-0.75){$>$}
\put(1.7,-0.7){\qbezier(0,0)(0.3,0)(1.2,0.7)}
\put(1.7,-0.7){\qbezier(0,0)(0.3,0)(2.3,0.7)}
\put(0.4,-0.15){$v_{1}$}
\put(1.7,0){\circle*{0.08}}
\put(1.6,0.1){$v_2$}
\put(2.9,0){\circle*{0.08}}
\put(2.8,0.1){$v_3$}
\put(4,0){\circle*{0.08}}
\put(3.9,0.1){$v_4$}
\put(4.2,0){\dots}
\put(2.6,0.5){\dots}
\put(2.6,-0.6){\dots}
\put(0.5,0.8){$e_1$}
\put(1.65,-0.85){$e_2$}
\end{picture}}
\vspace{2 cm}

In the example above $s^{-1}(r(e_i))^C=\{e_j:j<i\}$, which is finite, for each $i\in \N$.
\end{exemplo}

Notice that an ultragraph $\mathcal{G}$ with an infinite number of edges, and such that $s^{-1}(r(e))^C$ is finite for each edge $e$, has no sinks (that is, there is no $e\in \mathcal{G}^1$ such that each vertex $v\in r(e)$ is a sink). Moreover, each finite path $x\in \mathcal{G}^n$ is an element of the shift space $X$ associated to $\mathcal{G}$, so that, following Proposition \ref{p1}, 
$$X=\mathcal{G}^\infty\cup\bigcup\limits_{n=1}^\infty\mathcal{G}^n\cup \{\O\}.$$

\begin{remark} If $\mathcal{G}$ is an ultragraph with an infinite number of edges and such that $s^{-1}(r(e))^C$ is finite for each edge $e$, then $\mathcal{G}$ satisfies condition $(L)$ (that is, each loop has an exit). 
\end{remark}

Next we start to set up the ground for the definition of a partial action associated to an ultragraph $\mathcal{G}$.

Let $\F$ be the free group generated by $\mathcal{G}^1$ and let $0$ be the neutral element of $\F$. Notice that, for each $n\geq 1$, there is a copy of $\mathcal{G}^n=\{e_1...e_n:s(e_i)\in r(e_{i-1})\ \forall \ 2\leq i\leq n\}$  inside $\F$. Let $W\subseteq \F$ be the union of all these copies.

Define the following subsets of $X$:

\begin{itemize}
\item $X_0=X$,
\item $X_a=\{x\in X: x_1...x_{|a|}=a_1...a_{|a|}\}$, for each $a\in W$,
\item $X_{a^{-1}}=\{x\in X: s(x)\in r(a) \}\cup\{\O\}$, for all $a\in W$,
\item $X_{ab^{-1}}=\{x\in X:x=ay \text{ and }s(y)\in r(a)\cap r(b)\}\cup\{a\}$, for $ab^{-1}$ in the reduced form, with $a,b\in W$ and $r(a)\cap r(b)\neq \emptyset$.
\end{itemize}

\begin{proposicao} Let $\mathcal{G}$ be an ultragraph such that   $\mathcal{G}^1$ is infinite and $s^{-1}(r(e))^C$ is finite for each $e\in \mathcal{G}^1$. Then all the above sets $X_\alpha$ are clopen.
\end{proposicao}

\demo We prove the result for $X_{ab^{-1}}$, with $ab^{-1}$ in the reduced form, $a,b\in W$ and $r(a)\cap r(b)\neq \emptyset$ (notice that all other cases can be derived from this one, if we informally allow $a$ or $b$ to be the neutral element in $\F$).

To check that $X_{ab^{-1}}$ is closed, let $x^i$ be a sequence in $X_{ab^{-1}}$ converging to $x$. If $x=a$ we are done. If $x \neq a$ then $l(x)> l(a)$, so there exists $N>0$ such that, for all $n\geq N$ the $(l(a)+1)$ entry of $x^n$ is equal to the $(l(a)+1)$ entry of $x$. Hence $x\in X_{ab^{-1}}$.

Next we prove that each $X_{ab^{-1}}$ is open. Clearly, for each $x \in X_{ab^{-1}}$ different from $a$, any cylinder that contains $x$ is contained in $X_{ab^{-1}}$. So we just have to check that $a$ is an interior point of $X_{ab^{-1}}$. For this, notice that $F:=\left( s^{-1} \left(r(a) \cap r(b)\right) \right)^c$ is finite and then $Z(a,F)\subseteq X_{ab^{-1}}$ as desired. \fim

Associated to the sets $X_\alpha$ above we define the following maps:

\begin{itemize}
\item $\theta_0:X_0\rightarrow X_0$ by $\theta_0(x)=x$, for all $x\in X_0$,
\item $\theta_a:X_{a^{-1}}\rightarrow X_a$ by $\theta_a(\O)=a$, and $\theta_a(x)=ax$ for $x\neq \O $,
\item $\theta_{a^{-1}}:X_a\rightarrow X_{a^{-1}}$ by $\theta_{a^{-1}}(a)=\O$ and $\theta_{a^{-1}}(x)=x_{|a|+1}x_{|a|+2}...$ for $x\neq \O$,
\item $\theta_{ba^{-1}}:X_{ab^{-1}}\rightarrow X_{ba^{-1}}$ by $\theta_{ba^{-1}}(a)=b$ and $\theta_{ba^{-1}}
(x)=bx_{|a|+1}x_{|a|+2}...$ for $x\neq \O$.
\end{itemize}

Let $V\subseteq \F$ be the set $V=\{ab^{-1}:a,b\in W\cup\{0\}\}$. One can check that, for each $g\in V$, $\theta_g:X_{g^{-1}}\rightarrow X_g$ is a homeomorphism. Furthermore, $\theta_g(X_{g^{-1}}\cap X_h)\subseteq X_{gh}$ and $\theta_g(\theta_h(x))=\theta_{gh}(x)$, for each $x\in X_{(gh)^{-1}}\cap X_{h^{-1}}$ and $g,h\in V$.

Our aim is to define a partial action of $\F$ in the C*-algebra $C(X)$. With this in mind, for each $t\in V$, let $D_t=C(X_t)$ and define $\beta_t:D_{t^{-1}}\rightarrow D_t$ by $\beta_t(f)=f\circ \theta_{t^{-1}}$ and, for each $r\in \F\setminus V$, define $D_r$ as the null ideal of $C(X)$ and $\beta_r$ as the null map from $D_{r^{-1}}$ into $D_r$. We now have:

\begin{proposicao}\label{partialaction} $(\{D_t\}_{t\in \F}, \{\theta_t\}_{t\in \F})$ is a partial action of $\F$ in $C(X)$.
\end{proposicao}

\demo First note that for each $t\in V$, since $X_t$ is clopen and $\theta_t$ is a homeomorphism, $D_t$ is an ideal of $C(X)$ and $\beta_t:D_{t^{-1}}\rightarrow D_t$ is a C*-isomorphism. Moreover, $D_0=C(X)$ and $\beta_0=Id_{C(X)}$ by definition. The conditions $\beta_g(D_{g^{-1}}\cap D_h)\subseteq D_{gh}$ and $\beta_g(\beta_h(f))=\beta_{gh}(f)$, for each $f\in D_{(gh)^{-1}}\cap D_{h^{-1}}$, follow from the fact that $\theta_g(X_{g^{-1}}\cap X_h)\subseteq X_{gh}$ and $\theta_g(\theta_h(x))=\theta_{gh}(x)$, for each $x\in X_{h^{-1}}\cap X_{(gh)^{-1}}$ and $g,h\in V$. \fim

\begin{definicao} Let $\mathcal{G}$ be an ultragraph such that $\mathcal{G}^1$ is infinite and $s^{-1}(r(e))^C$ is finite for each $e\in \mathcal{G}^1$. The partial crossed product associated to $\mathcal{G}$ is the crossed product  $C(X)\rtimes_\beta\F$ arising from the partial action given in proposition \ref{partialaction}.
\end{definicao}

We prove in the next theorem that, under an additional assumption on $\mathcal{G}$, $C(X)\rtimes_\beta\F$ is isomorphic to the ultragraph C*-algebra $C^*(\mathcal{G})$. First we  prove the following lemma.

\begin{lema}\label{lemabeta} For each $t,a\in \F$ it holds that $\beta_t(1_{X_{t^{-1}}}1_{X_a})=1_{X_t}1_{X_{ta}}$
\end{lema}

\demo Given $y\in X$ and $U\subset X$, denote by $[y\in U]$ the binary function that is equal to 1 if $y\in U$ and 0 if $y\notin U$. Then, for each $x\in X$, we have that
$\beta_t(1_{X_{t^{-1}}}1_{X_a})(x)=(1_{X_{t^{-1}}}1_{X_a})(\theta_{t^{-1}}(x))=[\theta_{t^{-1}}(x)\in X_{t^{-1}}\cap X_a]=[x\in \theta_t(X_{t^{-1}}\cap X_a)]=[x\in X_t\cap X_{ta}]=1_{X_t}1_{X_{ta}}(x).$ \fim

\begin{remark}\label{setremark} From now on we assume that all ultragraphs have no sinks and, for each $v\in G^0$, $s^{-1}(v)$ is finite or $s^{-1}(v)^C$ is finite.
\end{remark}

The graph of Example \ref{ex1} satisfies the condition of the previous Remark. Below we give an example of an ultragraph that satisfies the conditions of Proposition \ref{partialaction}, but does not satisfy the condition of Remark \ref{setremark}.

\begin{exemplo} Let $G^0=\{v_i\}_{i\in \N}$, $\mathcal{G}^1=\{e_i:i\geq 2\}\cup \{e_{1j}:j\in \N\}$ and define $r(e_i)=\{v_j:j\geq i\}$ for $i\geq 2$; $r(e_{1j})=G^0$ for each $j\geq 1$; $s(e_i)=v_i$ for each $i\geq 2$ and $s(e_{1j})=v_1$ for each $j\geq 1$. We show a picture of this ultragraph below.

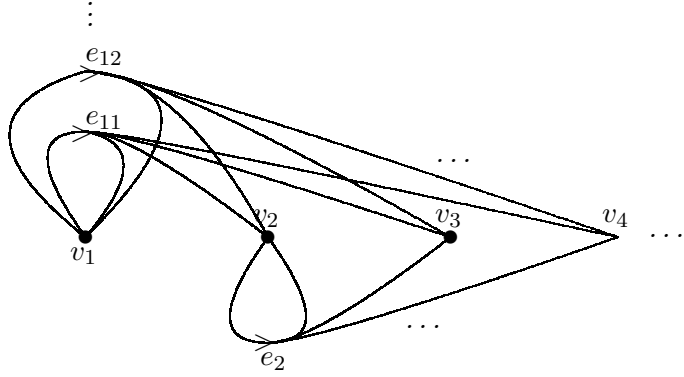
\centerline{
\setlength{\unitlength}{2cm}
\begin{picture}(4,1.7)
\put(0.5,0){\circle*{0.08}}
\put(0.5,0){\qbezier(0,0)(-0.5,0.7)(0,0.7)}
\put(0.5,0.7){\qbezier(0,0)(0.5,-0.1)(0,-0.7)}
\put(0.4,0.64){$>$}
\put(0.5,0.7){\qbezier(0,0)(0.3,0)(1.2,-0.7)}
\put(0.5,0.7){\qbezier(0,0)(0.3,0)(2.4,-0.7)}
\put(0.5,0.7){\qbezier(0,0)(0.3,0)(3.5,-0.7)}
\put(0.5,0){\qbezier(0,0)(-1,0.8)(0,1.1)}
\put(0.5,1.1){\qbezier(0,0)(1,-0.1)(0,-1.1)}
\put(0.45,1.035){$>$}
\put(0.5,1.17){$e_{12}$}
\put(0.5,1.1){\qbezier(0,0)(0.6,0)(1.2,-1.1)}
\put(0.5,1.1){\qbezier(0,0)(0.5,0)(2.4,-1.1)}
\put(0.5,1.1){\qbezier(0,0)(0.3,0)(3.5,-1.1)}
\put(1.7,0){\qbezier(0,0)(-0.5,-0.7)(0,-0.7)}
\put(1.7,0){\qbezier(0,0)(0.5,-0.7)(0,-0.7)}
\put(1.6,-0.75){$>$}
\put(1.7,-0.7){\qbezier(0,0)(0.3,0)(1.2,0.7)}
\put(1.7,-0.7){\qbezier(0,0)(0.3,0)(2.3,0.7)}
\put(0.4,-0.15){$v_{1}$}
\put(1.7,0){\circle*{0.08}}
\put(1.6,0.1){$v_2$}
\put(2.9,0){\circle*{0.08}}
\put(2.8,0.1){$v_3$}
\put(3.9,0.1){$v_4$}
\put(4.2,0){\dots}
\put(2.8,0.5){\dots}
\put(2.6,-0.6){\dots}
\put(0.5,1.4){\vdots}
\put(0.5,0.75){$e_{11}$}
\put(1.65,-0.85){$e_2$}
\end{picture}}
\vspace{2 cm}
\end{exemplo}

Notice that in the ultragraph above $s^{-1}(v_1)=\{e_{1j}:j\in \N\}$, which is infinite, and $s^{-1}(v_1)^C=\{e_i:i\geq 2\}$, which is also infinite. Furthermore, $s^{-1}(r(e))^C$ is finite for each edge $e$.

\begin{definicao}\label{setremark1} Let $\mathcal{G}$ be an ultragraph that satisfy all our standing hypothesis. Notice that in this case, for each $A\in \mathcal{G}^0$ it holds that either $s^{-1}(A)$ or $s^{-1}(A)^c$ is finite. Define, for each $A\in \mathcal{G}^0$,  $$X_A:=\{x\in X:s(x)\in A\},$$ if $s^{-1}(A)$ is finite and $$X_A:=\{x\in X:s(x)\in A\}\cup \{\O\}$$ if $s^{-1}(A)^C$ is finite. 
\end{definicao}

\begin{remark} Note that each $X_A$ is clopen, since $X_A=\bigcup\limits_{e\in s^{-1}(A)}X_e$, if $s^{-1}(A)$ is finite and $X_A=(\bigcup\limits_{e\in s^{-1}(A)^C}X_e)^C$, if $s^{-1}(A)^C$ is finite. Furthermore, for each $A,B\in \mathcal{G}^0$, we have that $X_A\cap X_B=X_{A\cap B}$ and $X_{A}\cup X_B=X_{A\cup B}$.
\end{remark}

\begin{teorema}\label{teoiso} Let $\mathcal{G}=(G^0,\mathcal{G}^1,r,s)$ be an ultragraph without sinks such that $\mathcal{G}^1$ is infinite, $s^{-1}(r(e))^C$ is finite for each $e\in \mathcal{G}^1$ and, for each $v\in G^0$, $s^{-1}(v)$ is finite or $s^{-1}(v)^C$ is finite. Then there exists a *-isomorphism $\Phi:C^*(\mathcal{G})\rightarrow C(X)\rtimes_\beta \F$ such that $\Phi(P_A)=1_{X_A}\delta_0$, for all $A\in \mathcal{G}^0$, and $\Phi(S_e)=1_{X_e}\delta_e$, for each $e\in \mathcal{G}^1$.
\end{teorema}

\demo We will define a *-homomorphism $\Phi:C^*(\mathcal{G})\rightarrow C(X)\rtimes_\beta \F$. Let, for each $e\in \mathcal{G}^1$, $\Phi(S_e):=1_{X_e}\delta_e$ and $\Phi(S_e^*):=1_{X_{e^{-1}}}\delta_{e^{-1}}$ and, for each $A\in \mathcal{G}^0$ define $\Phi(P_A):=1_{X_A}\delta_0$. To see that $\Phi$ extends to $C^*(\mathcal{G})$, we need to verify that $\Phi$ satisfies the relations of the universal C*-algebra $C^*(\mathcal{G})$. We do this below.

For each $A,B\in \mathcal{G}^0$, notice that $$\Phi(A)\Phi(B)=1_{X_A}\delta_01_{X_B}\delta_0=1_{X_A}1_{X_B}\delta_0=1_{X_{A}\cap X_B}\delta_0=1_{X_{A\cap B}}\delta_0=\Phi(P_{A\cap B})$$ and 
$$\Phi(P_{A\cup B})=1_{X_{A\cup B}}\delta_0=(1_{X_A}+1_{X_B}-1_{X_A\cap X_B})\delta_0=\Phi(P_A)+\Phi(P_B)-\Phi(P_{A\cap B}).$$

Now, for each $e\in \mathcal{G}^1$, notice that $$\Phi(S_e)^*\Phi(S_e)=1_{X_{e^{-1}}}\delta_{e^{-1}}1_{X_e}\delta_e=1_{X_{e^{-1}}}\delta_0=1_{X_{r(e)}}\delta_0=\Phi(P_{r(e)}).$$ 

Finally, since $\Phi(S_e)\Phi(S_e)^*=1_{X_e}\delta_0$, we have that $\{\Phi(S_e)\}_{e\in \mathcal{G}^1}$ have orthogonal ranges, $\Phi(S_e)\Phi(S_e)^*\leq \Phi(P_{s(e)})$ and $\Phi(P_v)=\sum\limits_{e\in s^{-1}(v)}\Phi(S_e)\Phi(S_e)^*$ if $0<|s^{-1}(v)|<\infty$. So, by the universal property of $C^*(\mathcal{G})$, we obtain a *-homomorphism $\Phi:C^*(\mathcal{G})\rightarrow C(X)\rtimes_\beta\F$.

To finish the proof we have to show that $\Phi$ is bijective. By \cite[17.11]{exel} we have that $\Phi(P_A)\neq 0$ and, since $\Phi$ is a *-homomorphism, $P_A\neq 0$. Since the ultragraph satisfies condition $(L)$, the Cuntz-Krieger Uniqueness theorem (see \cite[6.7]{Tom}), implies that $\Phi$ is injective.

It remains to show that $\Phi$ is surjective. Since, for each $t\in\F$,  $span\{1_{X_t}1_{X_a}:a\in W\cup\{0\}\}$ is dense in $D_t$, it is enough to show that $1_{X_t}1_{X_a}\delta_t\in Im(\Phi)$, where   $t=bc^{-1}$ belongs to $\F$, with $b,c\in W\cup\{0\}$, and $a\in W\cup\{0\}$. So, let $a,b,c\in W\cup\{0\}$. By applying the definition of $\Phi$ and Lemma \ref{lemabeta} it follows that $$\Phi(S_a)\Phi(S_a^*)\Phi(S_b)\Phi(S_c)^*=1_{X_a}\delta_0 1_{X_{bc^{-1}}}\delta_{bc^{-1}}=1_{X_a}1_{X_bc^{-1}}\delta_{bc^{-1}}$$ and hence $\Phi$ is surjective as desired. \fim

\begin{teorema} \label{isomorfismo} Let $\mathcal{G}$ and $\overline{\mathcal{G}}$ be two ultragraphs satisfying all the conditions of Theorem \ref{teoiso}. 
Let $X$ be the shift space associated to $\mathcal{G}$ and $Y$ be the shift space associated to $\overline{\mathcal{G}}$. If $X$ and $Y$ are conjugated then $C^*(\overline{\mathcal{G}})$ and $C^*(\mathcal{G})$ are isomorphic C*-algebras.
\end{teorema}

\demo By Theorem \ref{teoiso}, $C^*(\mathcal{G})$ is isomorphic to $C(X)\rtimes_\beta\F$.
So, it is enough to obtain a *-isomorphism $\Psi:C^*(\overline{\mathcal{G}})\rightarrow C(X)\rtimes_\beta\F$.
Let $\psi:Y\rightarrow X$ be a conjugacy. To obtain $\Psi$ we use the universality of $C^*(\overline{\mathcal{G}})$. Let us call the generators of $C^*(\overline{\mathcal{G}})$ by $\{p_A\}_{A\in \overline{\mathcal{G}}^0}$ and $\{s_e\}_{e\in \overline{\mathcal{G}}^1}$.
For each $A\in \overline{\mathcal{G}}^0$ define $\Psi(p_A)=1_{\psi(Y_A)}\delta_0$. 

Note that $$\begin{array}{ll} \Psi(p_A)\Psi(p_B) &= 1_{\psi(Y_A)}1_{\psi(Y_B)}\delta_0=1_{\psi(Y_A)\cap\psi(Y_B)}\delta_0=1_{\psi(Y_A\cap Y_B)}\delta_0 \\ & =1_{\psi(Y_{A\cap B})}\delta_0=\Psi(p_{A\cap B})\end{array}$$ 
and
$$\begin{array}{ll} \Psi(p_{A\cup B}) & =1_{\psi(Y_A\cup Y_B)}\delta_0=1_{\psi(Y_A)\cup \psi(Y_B)}\delta_0=(1_{\psi(Y_A)}+1_{\psi(Y_B)}-1_{\psi(Y_A)\cap \psi(Y_B)})\delta_0 \\ & = 
1_{\psi(Y_A)}\delta_0+1_{\psi(Y_B)}\delta_0-1_{\psi(Y_A\cap Y_B)}\delta_0=\Psi(p_A)+\Psi(p_B)-\Psi(p_{A\cap B}). \end{array}$$

We still need to define $\Psi$ on $S_e$. For this define, for each $e\in \overline{\mathcal{G}}^1$, $H_e:=\{a\in \mathcal{G}^1:X_a\cap \psi(Y_e)\neq\emptyset\}$ (which is a finite set, since $\psi(Y_e)$ is compact and $\{X_a\}_{a\in \mathcal{G}^1}$ is a disjoint open cover of $\psi(Y_e)$) and define $\Psi(s_e)=\sum\limits_{a\in H_e}1_{X_a}1_{\psi(Y_e)}\delta_a$. 

Next we check that $\{\Psi(p_A), \Psi(S_e) \}$ satisfy the remaining defining relations of C*($\overline{\mathcal{G}}^1$).

Note that 
$$\begin{array}{lll} \Psi(s_e)^*\Psi(s_e)= & \sum\limits_{a\in H_e}\beta_{a^{-1}}(1_{X_a}1_{\psi(Y_e)})\delta_{a^{-1}}\sum\limits_{b\in H_e}1_{X_b}1_{\psi(Y_e)}\delta_b \\
& =\sum\limits_{a,b\in H_e}\beta_{a^{-1}}(\beta_a(\beta_a^{-1}(1_{X_a}1_{\psi(Y_e)}))1_{X_b}1_{\psi(Y_e)})\delta_{a^{-1}b} \\
& =\sum\limits_{a\in H_e}\beta_{a^{-1}}(1_{X_a}1_{\psi(Y_e)})\delta_0\end{array}.$$

The equality $\Psi(s_e)^*\Psi(s_e)=\Psi(p_{r(e)})$ now comes from the following claim:

{\it Claim 1: $\sum\limits_{a\in H_e}\beta_{a^{-1}}(1_{X_a}1_{\psi(Y_e)})=1_{\psi(Y_{r(e)})}$}

Let $x\in X$. Then
$$\sum\limits_{a\in H_e}\beta_{a^{-1}}(1_{X_a}1_{\psi(Y_e)})(x)=\sum\limits_{a\in H_e}1_{X_{a^{-1}}}(x)(1_{X_a}1_{\psi(Y_e)})(\theta_a(x))=$$
$$=\sum\limits_{a\in H_e}[x\in X_{a^{-1}}][\theta_a(x)\in \psi(Y_e)]=\sum\limits_{a\in H_e}[x\in X_{a^{-1}}][ax\in \psi(Y_e)]=$$
$$=\sum\limits_{a\in H_e}[x\in X_{a^{-1}}][\psi^{-1}(ax)\in Y_e].$$ Since $\psi^{-1}$ is a conjugacy then $\psi^{-1}(ax)=z_a\psi^{-1}(x)$ for some $z_a\in \overline{\mathcal{G}}^1$ and so $\psi^{-1}(ax)\in Y_e$ if, and only if, $z_a=e$. Then 
$\sum\limits_{a\in H_e}[x\in X_{a^{-1}}][\psi^{-1}(ax)\in Y_e]=0$ or $\sum\limits_{a\in H_e}[x\in X_{a^{-1}}][\psi^{-1}(ax)\in Y_e]=1.$ 
We want to show that $$\sum\limits_{a\in H_e}[x\in X_{a^{-1}}][\psi^{-1}(ax)\in Y_e]=[x\in \psi(Y_{r(e)})],$$ for each $x\in X$. 

Fix an $x\in X$. Suppose that $\sum\limits_{a\in H_e}[x\in X_{a^{-1}}][\psi^{-1}(ax)\in Y_e]=1$. If $x=\O$ then $[x\in \psi(Y_{r(e)})]=1$, since $\O\in Y_{r(e)}$. So, let $x\neq \O$. Then $\psi^{-1}(ax)=z_a\psi^{-1}(x)= e\psi^{-1}(x)$, 
from where it follows that $s(\psi^{-1}(x))\in r(e)$ and hence $[x\in \psi(Y_{r(e)})]=1$. 

Now, let $x\in X$ be such that $[x\in \psi(Y_{r(e)})]=1$. If $x=\O$ then one of the elements of $H_e$ is $a=\psi(e)$ and, for this $a$, we have that $[\O\in X_{a^{-1}}][\psi^{-1}(a\O)\in Y_e]=1$. Furthermore, for $x=\O$, $\psi^{-1}(ax)$ has lenght one for every $a\in H_e$, so that $[\psi^{-1}(ax) \in Y_e]=0$ for all $a\neq \psi(e)$. So, we may suppose $x\neq \O$. Since $[x\in \psi(Y_{r(e)})]=1$ then $\psi^{-1}(x)\in Y_{r(e)}$ and hence $e\psi^{-1}(x)\in Y_e$, from where we obtain that $\psi(e\psi^{-1}(x))\in \psi(Y_e)$. Since $\psi(Y_e)\subseteq \bigcup\limits_{a\in H_e}X_a$ then $\psi(e\psi^{-1}(x))\in X_{a_e}$, for some $a_e\in H_e$, what implies that $\psi(e\psi^{-1}(x))=a_e x$ for some $x\in X_{a_e^{-1}}$.
Moreover, $\psi^{-1}(a_ex)=e\psi^{-1}(x)\in Y_e$, so that $[x\in X_{a_e^{-1}}][\psi^{-1}(a_ex)\in Y_e]=1$ and hence $\sum\limits_{a\in H_e}[x\in X_{a^{-1}}][\psi^{-1}(ax)\in Y_e]=1$. Therefore the claim is proved.

To verify that $\{\Psi(p_A), \Psi(S_e) \}$ satisfy relations 3 and 4 defining C*($\overline{\mathcal{G}}^1$) we need the following auxiliary result:

{\it Claim 2: For each $e\in \overline{\mathcal{G}}^1$, and for each $w\in \bigcup\limits_{n=1}^\infty\overline{\mathcal{G}}^n$ such that $s(e)\in r(e)$, it holds that $\Psi(s_e)\Psi(s_e)^*=1_{\psi(Y_e)}\delta_0$ and $\Psi(s_e)1_{\psi(Y_w)}\Psi(s_e)^*=1_{\psi(Y_{ew})}\delta_0$}

First we compute $\Psi(s_e)1_{\psi(Y_w)}\Psi(s_e)^*$:
$$\begin{array}{ll} Psi(s_e)\Psi(s_e)^* & =\sum\limits_{a\in H_e}1_{X_a}1_{\psi(Y_e)}\delta_a\sum\limits_{b\in H_e}1_{\psi(Y_w)}\beta_{b^{-1}}(1_{X_b}1_{\psi(Y_e)})\delta_{b^{-1}} \\
& =\sum\limits_{a,b\in H_e}\beta_a(\beta_{a^{-1}}(1_{X_a}1_{\psi(Y_e)})1_{\psi(Y_w)}\beta_{b^{-1}}(1_{X_b}1_{\psi(Y_e)}))\delta_{ab^{-1}}.\end{array}$$ Note that for each $x=ay\in X_{ab^{-1}}$,
$$\beta_a(\beta_{a^{-1}}(1_{X_a}1_{\psi(Y_e)})1_{\psi(Y_w)}\beta_{b^{-1}}(1_{X_b}1_{\psi(Y_e)}))(x)=$$ $$=1_{X_a}(ay)1_{\psi(Y_e)}(ay)1_{\psi(Y_w)}(y)1_{X_b}(by)1_{\psi(Y_e)}(by)=$$
$$=1_{\psi(Y_e)}(ay)1_{\psi(Y_w)}(y)1_{\psi(Y_e)}(by)=$$ $$=[\psi^{-1}(ay)\in Y_e][\psi^{-1}(y)\in Y_w][\psi^{-1}(by)\in Y_e].$$ 
For $a\neq b$, if $\psi^{-1}(ay)\in Y_e$ then $\psi^{-1}(ay)=e\psi^{-1}(y)$ and so $\psi^{-1}(by)\notin Y_e$, because otherwise $\psi^{-1}(by)=e\psi^{-1}(y)=\psi^{-1}(ay)$, which is impossible since $\psi^{-1}$ is injective. 

Hence, for $a\neq b$ we have that $$\beta_a(\beta_{a^{-1}}(1_{X_a}1_{\psi(Y_e)})1_{\psi(Y_w)}\beta_{b^{-1}}(1_{X_b}1_{\psi(Y_e)}))(x)=0.$$ On the other hand, if $a=b$ then $$\beta_a(\beta_{a^{-1}}(1_{X_a}1_{\psi(Y_e)})1_{\psi(Y_w)}\beta_{b^{-1}}(1_{X_b}1_{\psi(Y_e)}))(x)= 1_{\psi(Y_e)}(x)1_{X_a}(x)1_{\psi(Y_w)}(y).$$ 

Now, notice that for $x=ay\in X_a$ we have that $$1_{\psi(Y_e)}(x)1_{X_a}(x)1_{\psi(Y_w)}(y)=[ay\in \psi(Y_e)][y\in \psi(Y_w)]=$$
$$=[\psi^{-1}(ay)\in Y_e][\psi^{-1}(y)\in Y_w]=[\psi^{-1}(ay)\in Y_{ew}]=$$ $$=[ay\in \psi(Y_{ew})]=1_{\psi(Y_{ew})}(x)$$ and hence, for $a=b$, we have that
$$\beta_a(\beta_{a^{-1}}(1_{X_a}1_{\psi(Y_e)})1_{\psi(Y_w)}\beta_{b^{-1}}(1_{X_b}1_{\psi(Y_e)}))(x)=1_{X_a}1_{\psi(Y_{ew})}(x),$$ for all $x\in X$. Therefore, 
$$\Psi(s_e)1_{\psi(Y_w)}\Psi(s_e)^*=\sum\limits_{a\in H_e}1_{X_a}1_{\psi(Y_{ew})}\delta_0=1_{\psi(Y_{ew})}\delta_0,$$ where the last equality follows from the fact that $\psi(Y_{ew})\subseteq\psi(Y_e)\subseteq\bigcup\limits_{a\in H_e}X_a$. 

The equality $\Phi(s_e)\Phi(s_e)^*=1_{\psi(Y_e)}\delta_0$ follows analogously, by replacing $1_{\psi(Y_w)}$ by $1_X$. This finishes the proof of Claim 2.

Now, notice that Claim 2 readily implies that $\Psi(s_e)\Psi(s_e)^*\leq \Psi(p_{s(e)})$, since $Y_e\subseteq Y_{s(e)}$. Moreover, also by Claim 2, $\{\Psi(s_e)\}_{e\in \overline{\mathcal{G}}^1}$ is a family of partial isometries with orthogonal ranges.

Finally, for $v\in \overline{\mathcal{G}}^0$ such that $0<|s^{-1}(v)|<\infty$, since $Y_v=\bigcup\limits_{e\in s^{-1}(v)}Y_e$, we have that $$\Psi(p_v)=1_{\psi(Y_v)}\delta_0=\sum\limits_{e\in s^{-1}(v)}{1_{\psi(Y_e)}}\delta_0=\sum\limits_{e\in s^{-1}(v)}\Psi(s_e)\Psi(s_e)^*.$$

By the universality of $C^*(\mathcal{G})$, we obtain a *-homomorphism, which we also call $\Psi$, from  $C^*(\mathcal{G})$ to $C(X)\rtimes_\beta\F$. By \cite[17.11]{exel}, $\Psi(p_A)=1_{\psi(Y_A)}\delta_0\neq 0$ for each $A\in \mathcal{G}^0$. Since $\Psi$ is a *-homomorphism then $p_A\neq 0$ for each $A\in \mathcal{G}^0$. So, it follows from \cite[6.7]{Tom} that $\Psi$ is *-isomorphism from $C^*(\overline{\mathcal{G}})$ to the sub-algebra of $C(X)\rtimes_\beta\F$ generated by $\{\Psi (p_A)\}_{A\in \overline{\mathcal{G}}}$ and $\{\Psi(s_e)\}_{e\in \overline{\mathcal{G}}^1}$.

It remains to prove that $\Psi$ is surjective. 

Note that $span\{1_Y,1_{Y_w}:w\in \bigcup\limits_{n=1}^\infty\overline{\mathcal{G}}^n\}$ is dense in $C(Y)$ and, since $\psi:Y\rightarrow X$ is a homeomorphism, we have that $span\{1_X,1_{\psi(Y_w)}:w\in \bigcup\limits_{n=1}^\infty\overline{\mathcal{G}}^n\}$ is dense in $C(X)$.
Let $w=e_1...e_m\in \bigcup\limits_{n=1}^\infty\overline{\mathcal{G}}^n$. By induction over the length of $w$, together with Claim 2, we obtain that $\Psi(s_{e_1})...\Psi(s_{e_m})\Psi(s_{e_m})^*...\Psi(s_{e_1})^*=1_{\psi(Y_w)}\delta_0$. Moreover, $\Psi(1)=\Psi(p_{\overline{G}^0})=1_{\psi(Y_{\overline{G}^0})}\delta_0=1_{\psi(Y)}\delta_0=1_X\delta_0$. 

By the density of $span\{1_X,1_{\psi(Y_w)}:w\in \bigcup\limits_{n=1}^\infty\overline{\mathcal{G}}^n\}$ in $C(X)$ it follows that $C(X)\delta_0\subseteq Im(\Psi)$. In particular, $1_{X_a}\delta_0 \in Im(\Psi)$ for each $a\in \mathcal{G}^1$.

Next, fix a $b\in \mathcal{G}^1$. Note that $\{\psi(Y_e)\}_{e\in \overline{\mathcal{G}}^1}$ is an open disjoint cover of the compact set $X_b$. Let $H:=\{e\in \overline{\mathcal{G}}^1:\psi(Y_e)\cap X_b\neq\emptyset\}$, which is finite since $X_b$ is compact.
Then $\sum\limits_{e\in H}1_{X_b}\delta_0\Psi(s_e)\in Im(\Psi)$ and 
$$\sum\limits_{e\in H}1_{X_b}\delta_0\Psi(s_e)=\sum\limits_{e\in H}\sum\limits_{a\in H_e}1_{X_b}1_{X_a}1_{\psi(Y_e)}\delta_a=\sum\limits_{e\in H}1_{X_b}1_{\psi(Y_e)}\delta_b=1_{X_b}\delta_b,$$ where the last equality is true since $X_b\subseteq \bigcup\limits_{e\in H}\psi(Y_e)$.

So, it follows that $\{1_{X_a}\delta_a:a\in \mathcal{G} ^1\}\subseteq Im(\Psi)$ and, since $C(X)\delta_0\subseteq Im(\Psi)$, in particular $\{1_{X_A}\delta_0:A\in \mathcal{G}^0\}\subseteq Im(\Psi)$. By Theorem \ref{teoiso}, the elements $\{1_{X_a}\delta_a\}_{a\in \mathcal{G}^1}$ and $\{1_{X_A}\delta_0\}_{A\in \mathcal{G}^0}$ generate $C(X)\rtimes_\beta\F$ and hence $\Psi$ is surjective. \fim


\section{Shifts of finite type and ultragraphs.}

In this section we produce two shifts of finite type that are not eventually finite periodic conjugate to an edge shift of a graph. We also show that a graph edge shift of finite type is always a full shift, while this is not true for ultragraph edge shifts of finite type. Before we proceed we recall the reader the methods used to compute the K-theory of ultragraph C*-algebras, as described in \cite{KMST}. 

Let $\mathcal{G}=(G^0, \mathcal{G}^1, r,s)$ be an ultragraph. We denote by $Z_{\mathcal{G}}$ the (algebraic) subalgebra of $l^{\infty}(G^0, \Z)$ generated by $\{\delta_v: v\in G^0\}\cup \{\chi_{r(e)}: e \in \mathcal{G}^1 \}$. In \cite{KMST} it is proved that $Z_{\mathcal{G}}= \{ \displaystyle \sum_{k=1}^{n} z_k \chi_{A_k}: \ n\in \N, z_k \in \Z, A_k \in \mathcal{G}^0\}$. A vertex $v\in G^0$ is said to be regular if $0<|s^{-1}(v)|<\infty$ and the set of all regular vertices is denoted by $G^0_{rg}$. 

\begin{teorema}\cite[Theorem~5.4]{KMST}\label{Ktheory} Let $\mathcal{G}=(G^0, \mathcal{G}^1, r,s)$ be an ultragraph. We define $\partial:\Z^{G^0_{rg}} \rightarrow Z_{\mathcal{G}}$ by $\partial (\delta_v) = \delta_v - \sum_{e\in s^{-1}(v)} \chi_{r(e)}$ for $v \in G^0_{rg}$. Then we have $K_0 \left(C^*(\mathcal{G})  \right) = \text{coker}(\partial)$ and $K_1 \left(C^*(\mathcal{G})  \right) = \text{ker}(\partial)$.
\end{teorema}

Notice that for the full shift, which can be represented by the edge shift of a graph with one vertex and countable number of loops, the $K_0$ group of the associated C*-algebra is $\Z$, since the map $\partial$ of definition \ref{Ktheory} is the zero map and there are no regular vertices in this graph.

Let $\mathcal{G}$ be the ultragraph of example \ref{exemplo1}, that is $\mathcal{G}$ is the ultragraph with a countable number of vertices, say $G^0=\{v_i\}_{i\in \N}$, and edges such that $s(e_i)= v_i$ for all $i$, $r(e_1)=\{v_3, v_4, v_5, \ldots \}$ and $r(e_j)= G^0$ for all $j\neq 1$.  Our goal is to show that the $K_0$ group of $C^*(\mathcal{G})$ is different from $\Z$. We will show, using the theorem above, that this group has torsion. 

\begin{proposicao}\label{torsao} $K_0(C^*(\mathcal{G})) = \Z \oplus \Z_2$.
\end{proposicao}

\demo To compute $K_0$ we need to compute the map $\partial$ and its cokernel. Notice that in $\mathcal{G}$ all vertices are regular, since they emit only one edge. 

Now, $\partial(\delta_{v_1}) = \delta_{v_1} - \chi_{G^0\setminus\{v_1,v_2\}}$ and, for each vertex $v_j$ with $j\neq 1$, we have that $\partial(\delta_{v_j}) = \delta_{v_j} - \chi_{G^0}$. Therefore the function $\chi_{G^0}$ is equivalent to $\delta_{v_j}$, which is equivalent to $\delta_{v_k}$ for all $j, k \neq 1$, and $\delta_{v_1}$ is equivalent to $\chi_{G^0\setminus\{v_1,v_2\}}$. Furthermore, $$\partial(\delta_{v_1} - \delta_{v_2}) = \delta_{v_1} - \delta_{v_2} - \chi_{G^0 \setminus\{v_1,v_2\}} + \chi_{G^0} = \delta_{v_1} - \delta_{v_2} + \delta_{v_1} + \delta_{v_2} = 2 \delta_{v_1},$$
so that $2 \delta_{v_1} \in \text{Im}(\partial)$ and hence $2 \delta_{v_1} = 0$ in $\text{coker}(\partial)$. 

From the observations above we see that $\text{coker}(\partial)$ is generated by the class of $\delta_{v_1}$ and $\delta_{v_2}$. To finish the proof all we need to do is show that the classes $\delta_{v_1}$ and $n \delta_{v_2}$, $n\in \Z$, are non-zero, that is, we need to show that $\delta_{v_1}$ and $n\delta_{v_2}$, $n\in \Z$, are not in the image of $\partial$.

Suppose $\delta_{v_1} \in \text{Im}(\partial)$. Then there exists $c_i\in \Z$, which eventually vanish, such that \begin{equation}\label{eqK} \delta_{v_1} = \sum c_i \partial(\delta_{v_i}) =  c_1 \delta_{v_1} -c_1 \chi_{G^0\setminus\{v_1,v_2\}} + \sum_{i\neq 1} c_i \delta_{v_i} - c_i \chi_{G^0}. \end{equation}
Evaluating the above equation at $v_j$ we obtain that $$\begin{array}{lll} 1 = c_1 - \sum_{i\neq 1} c_i, \\ 
0 = c_2 - \sum_{i\neq 1} c_i, \\
0 = -c_1 + c_j - \sum_{i\neq 1} c_i, \text{for } j\neq 2. 
\end{array}$$
The last equation above implies that $c_j=c_k$ for all $j\neq 3$ and, since they vanish eventually, this implies that $c_j=0$ for all $j\geq 3$. Also this equation now implies that $c_1+c_2=0$ and the first equation gives us that $1=c_1-c_2$. Therefore $c_2 = -\frac{1}{2}$ and $c_1=\frac{1}{2}$, a contradiction (since all $c_i$ should be integers).

To prove that for every for non-zero integer $n$ the function $n\delta_{v_2}$ is not in the image of $\partial$ one proceeds analogously, evaluating equation \ref{eqK} (with the left side replaced with $n\delta_{v_2}$) on the vertices $\{v_i \}$ and deriving a contradiction. \fim

\begin{exemplo}
\label{exemplo2} Let $\mathcal{G}_1$ be the ultragraph with a countable number of vertices, say $G^0=\{v_i\}_{i\in \N}$, and edges $e_0, e_1, e_2, \ldots$ such that $s(e_0)=v_1$, $s(e_i)= v_i$ for all $i>0$, $r(e_0)=G^0\setminus\{v_2\}$, $r(e_1)=r(e_2)=G^0\setminus\{v_1,v_2\}$ and $r(e_j)= G^0$ for all $j>2$. 
\end{exemplo}

\begin{proposicao} The $K_0$ group of $(C^*(\mathcal{G}_1)) $ has torsion.
\end{proposicao}

\demo As before we need to compute the map $\partial$. Notice that $\partial(\delta_{v_1}) = \delta_{v_1} - \chi_{G^0\setminus\{v_2\}}- \chi_{G^0\setminus\{v_1,v_2\}} = -2 \chi_{G^0\setminus\{v_1,v_2\}}$, $\partial(\delta_{v_2}) = \delta_{v_2} - \chi_{G^0\setminus\{v_1,v_2\}} $ and, for each vertex $v_j$ with $j>2$, $\partial(\delta_{v_j}) = \delta_{v_j} - \chi_{G^0}$. 

From the above it is clear that the class of $\delta_{v_2}$ is the same as the class of $\chi_{G^0\setminus\{v_1,v_2\}}$ in the cokernel of $\partial$ and that this class has torsion. All is left to do is to show that $\delta_{v_2}$ is not in the image of $\partial$. So suppose that $\delta_{v_2} \in \text{Im}(\partial)$. Then there exists $c_i\in \Z$, which eventually vanish, such that 
\begin{equation}\label{eq2} \delta_{v_2} = c_1 (-2 \chi_{G^0\setminus\{v_1,v_2\}}) + c_2(\delta_{v_2} - \chi_{G^0\setminus\{v_1,v_2\}}) + \sum_{i> 2} c_i \delta_{v_i} - c_i \chi_{G^0}. \end{equation}
Evaluating the above equation at $v_j$, $j>2$, we obtain that all $c_k=0$ for $k>2$. Also we get that $-2c_1-c_2=0$ and from equation \ref{eq2} evaluated at $v_2$ we have that $c_2=1$. These combined imply that $c_1=-\frac{1}{2}$, which is not an integer. Therefore $\delta_{v_2}$ is not in the image of $\partial$ as desired. \fim

We now have all the ingredients to prove our main theorem:

\begin{teorema}\label{provafracaconjectura} Let $A$ be an infinite alphabet and $\Sigma_A$ be the one sided full shift as in \cite{OTW}. Then there exists subshifts of finite type of $\Sigma_A$ that are not conjugate, via an eventually finite periodic conjugacy, to an edge shift of a graph. In particular, the shifts of finite type given as the edge shifts of the ultragraphs in examples \ref{exemplo1} and \ref{exemplo2} are not eventually finite periodic conjugate to an edge shift of a graph.
\end{teorema}

\demo Let $A=\{a_i\}$ be a countable alphabet, $F=\{a_1a_1, a_1a_2\}$ and consider the shift of finite type $X_F$ as described in section 2.1. It is straightforward to check that $X_F$ is conjugate to the edge shift of the ultragraph $\mathcal{G}$ given in example \ref{exemplo1}. 

Suppose that $X_F$ is eventually finite periodic conjugate to an edge shift of a graph. Then, by proposition \ref{fullshift}, the edge shift of the graph is a full shift and therefore conjugate to the edge shift of a graph with one vertex and countable number of loops. By theorem \ref{isomorfismo}, the C*-algebras associated to the ultragraph of example \ref{exemplo1} and the graph describing the full shift are isomorphic and hence have the same K-theory. But $K_0(C^*(\mathcal{G}))$ has torsion, as showed in Proposition \ref{torsao}, while $K_0$ of the graph associated to the full shift is $\Z$. We obtain therefore a contradiction and hence $X_F$ is not eventually finite periodic conjugate to the edge shift of a graph. 

The proof that the shift arising as the edge shift of the ultragraph given in example \ref{exemplo2} is not eventually finite periodic conjugate to the edge shift of a graph is analogous to what is done above. \fim

Next we show that a graph edge shift of finite type is always a full shift.

\begin{proposicao} Let $X_E$ be an edge shift of a graph $E$, where $E^1$ is infinite. Suppose $X_E$ is a shift of finite type. Then $X_E$ is the full shift $\Sigma_{E^1}$.
\end{proposicao}

\demo Recall that $X_E$ is equal to $X_F$, where $F:=\{ef: r(e)\neq s(f)$. Also, since $X$ is a SFT, there exists a finite set of forbidden words $G$ such that $X=X_G$. Since $X_E$ is a SFT, each forbidden word $ef \in F$ must be a subblock of some forbidden word in $G$ (otherwise there would be infinitely many $g_i\in E^1$ such that $efg_i \in X_G$ and hence $ef \in X_G$). So $F$ is also a finite set. 

Now, let $e,f\in E^1$. Since $E^1$ is infinite and $F$ is finite, there exists edges $e_1$ and $e_2$ such that $r(e)=s(e_1)=r(f)$ and $s(e)=r(e_2)=s(f)$. So, all the edges have the same range and the same source. Again, since $F$ is finite, there exists at least one loop in $E$ and hence all the edges are loops based on the same vertex. Hence $X_E=\Sigma_{E^1}$. \fim

We finalize the paper showing that the result above is not true for ultragraph edge shifts of finite type. 

\begin{proposicao} There exists an ultragraph edge shift of finite type which is not conjugated to any full shift.
\end{proposicao}

\demo Let $X_\mathcal{G}$ be the edge shift of the ultragraph $\mathcal{G}$ of Example \ref{exemplo1}. Note that $X_\mathcal{G}$ is a shift of finite type (the set of forbidden words is $\{e_1e_1,e_1e_2\}$). Suppose that $X_\mathcal{G}$ is conjugated to a full shift $\Sigma_A$. 
As before, we can see $\Sigma_A$ as the edge shift space associated to the graph $E$ consisting of a bouquet of loops with same range. By Theorem \ref{teoiso}, $C^*(\mathcal{G})$ and $C^*(E)$ are isomorphic C*-algebras, but this is impossible, since $K_0(C^*(E))=Z$ and, by Proposition \ref{torsao}, $K_0(C^*(\mathcal{G})=\Z\oplus \Z_2$. Hence $X_G$ is not conjugated to a full shift. \fim

\end{document}